\documentclass[12pt]{amsart}
\usepackage{amsmath}
\usepackage{amssymb}
\usepackage{latexsym}

\newcommand{\Zint}{\mathbb {Z}}    
     
\newcommand{\Rea}{\mathbb {R}}      % Real number field
\newcommand{\Cplx}{\mathbb {C}}     % Complex  number field

\newcommand{\halmos}{\rule{5pt}{5pt}}

\numberwithin{equation}{section}

\newtheorem{prop}{\bf Proposition}[section]

\begin{document}

\title[Heun equation and Painlev\'e equation]
{Heun equation and Painlev\'e equation}
\author{Kouichi Takemura}
\address{Department of Mathematical Sciences, Yokohama City University, 22-2 Seto, Kanazawa-ku, Yokohama 236-0027, Japan.}
\email{takemura@yokohama-cu.ac.jp}
\thanks{Partially supported by a Grant-in-Aid for Scientific Research (No. 15740108) from the Japan Society for the Promotion of Science.}

\subjclass{33E10,34M55}

\begin{abstract}
We relate two parameter solutions of the sixth Painlev\'e equation and finite-gap solutions of the Heun equation by considering monodromy on a certain class of Fuchsian differential equations. In the appendix, we present formulae on differentials of elliptic modular functions, and obtain the ellitic form of the sixth Painlev\'e equation directly.
\end{abstract}

\maketitle

\section{Introduction}

In this paper we make a study on two differential equations. One is the Heun equation, and the other is the sixth Painlev\'e equation.

Heun's differential equation (or the Heun equation) is a differential equation given by
\begin{equation}
\left( \! \left(\frac{d}{dw}\right) ^2 \! + \left( \frac{\gamma}{w}+\frac{\delta }{w-1}+\frac{\epsilon}{w-t}\right) \frac{d}{dw} +\frac{\alpha \beta w -q}{w(w-1)(w-t)} \right)\tilde{f}(w)=0
\label{Heun}
\end{equation}
with the condition 
\begin{equation}
\gamma +\delta +\epsilon =\alpha +\beta +1.
\label{Heuncond}
\end{equation}
The Heun equation is the standard canonical form of a Fuchsian equation with four singularities.
It is well known that the Fuchsian equation with three singularities is the hypergeometric differential equation.

In the 1980's, Treibich and Verdier \cite{TV} found that the Heun equation is related with the theory of the finite-gap potential, and several others have produced more precise statements and concerned results on this subject.
Namely, integral representations of solutions, global monodromy in terms of hyperelliptic integrals and the Hermite-Krichever Ansatz for the case $\gamma, \delta, \epsilon , \alpha -\beta \in \Zint +\frac{1}{2}$ are investigated (see \cite{BE,GW,Smi,Tak1,Tak2,Tak3,Tak4} etc.).

The sixth Painlev\'e equation is a non-linear ordinary differential equation written as
\begin{align}
\frac{d^2\lambda }{dt^2} = & \frac{1}{2} \left( \frac{1}{\lambda }+\frac{1}{\lambda -1}+\frac{1}{\lambda -t} \right) \left( \frac{d\lambda }{dt} \right) ^2 -\left( \frac {1}{t} +\frac {1}{t-1} +\frac {1}{\lambda -t} \right)\frac{d\lambda }{dt} \label{eq:P6eqnintr} \\
& +\frac{\lambda (\lambda -1)(\lambda -t)}{t^2(t-1)^2}\left\{ \frac{\kappa _{\infty}^2}{2} -\frac{\kappa _{0}^2}{2}\frac{t}{\lambda ^2} +\frac{\kappa _{1}^2}{2}\frac{(t-1)}{(\lambda -1)^2} +\frac{(1-\kappa _{t}^2)}{2}\frac{t(t-1)}{(\lambda -t)^2} \right\}. \nonumber
\end{align}
A remarkable property of this differential equation is that its solutions do not have movable singularities other than poles.
Although generic solutions of the sixth Painlev\'e equation are trancedental, it may have classical solutions for special cases.
If $\kappa _0  = \kappa _1 = \kappa _t = \kappa _{\infty} =0$, then Eq.(\ref{eq:P6eqnintr}) has two parameter solutions called Picard's solution \cite{Oka}, and if $\kappa _0  = \kappa _1 = \kappa _t = \kappa _{\infty} =1/2$, then Eq.(\ref{eq:P6eqnintr}) has two parameter solutions called Hitchin's solution \cite{Hit}.

In this paper we investigate a family of solutions to the sixth Painlev\'e equation including Hitchin's solutions by applying Hermite-Krichever Ansatz which is used to study the Heun equation in \cite{Tak4}.
More precisely, we develop the Hermite-Krichever Ansatz for a certain class of Fuchsian differential equations which include the linear differential equation that produce the sixth Painlev\'e equation by monodromy preserving deformation.
By considering monodromy preserving deformation for the solutions to the linear differential equation, we obtain solutions to the sixth Painlev\'e equation including Hitchin's solutions. 

This paper is organized as follows. 
In section \ref{sec:HK}, we obtain integral representations of solutions to a certain class of Fuchsian differential equations and rewrite them to the form of the Hermite-Krichever Ansatz.
In section \ref{sec:Heun}, we apply the results in section \ref{sec:HK} for the Heun equation. 
In section \ref{sec:P6}, we show that the solutions to the linear differential equations considered in section \ref{sec:HK} produce two parameter solutions to the sixth Painlev\'e equation by monodromy preserving deformation. Some explicit solutions that include Hitchin's solution are displayed.
In section \ref{sec:rmk}, we give concluding remarks and present an open problem.
In the appendix, we present formulae on differentials of elliptic modular functions, and obtain the ellitic form of the sixth Painlev\'e equation directly.

\section{Fuchsian differential equation and Hermite-Krichever Ansatz} \label{sec:HK}
In this section, we consider differential equations which have additional apparent singularities to the Heun equation. More precisely, we consider the equation
\begin{align}
& \left\{ \frac{d^2}{dw^2}+\left( \frac{\frac{1}{2}-l_1}{w}+  \frac{\frac{1}{2}-l_2}{w-1}+  \frac{\frac{1}{2}-l_3}{w-t}+ \sum _{{i'}=1}^M \frac{-r_{i'}}{w-\tilde{b}_{i'}} \right) \frac{d}{dw} \right. \label{Feq} \\
& \left.  +\frac{(\sum _{i=0}^3 l_i + \sum _{{i'}=1}^M r_{i'})(-1-l_0 +\sum _{i=1}^3 l_i + \sum _{{i'}=1}^M r_{i'})w+\tilde{p}+ \sum_{{i'}=1}^M \frac{\tilde{o}_{i'}}{w-\tilde{b}_{i'}}}{4w(w-1)(w-t)}\right\}\tilde{f}(w) =0 , \nonumber
\end{align}
for the case $l_i \in \Zint _{\geq 0}$ $(0 \leq i \leq 3)$, $r_{i'} \in \Zint _{>0}$ $(1 \leq i' \leq M)$ and the regular singular points $\tilde{b}_{i'}$ $(1 \leq i' \leq M)$ are apparent.
Here, a regular singular point $x=a$ of a linear differential equation of order two is said to be apparent, if and only if the differential equation does not have a logarithmic solution at $x=a$ and the exponents at $x=a$ are integers.

Let $\wp (x)$ be the Weierstrass $\wp$-function with periods $(2\omega_1, 2\omega _3)$.  We set $\omega _0 =0$, $\omega _1 =1/2$ $\omega _3 =\tau /2$, $\omega _2=-\omega _1-\omega _3$ and $e_i=\wp (\omega _i)$ $(i=1,2,3)$.
It is known that, if $t \neq 0,1,\infty$, then there exists a value $\tau \in \Rea + \sqrt{-1} \Rea _{>0}$ such that $t=(e_3-e_1)/(e_2 -e_1)$.
By a certain transformation, Eq.(\ref{Feq}) is rewritten in terms of elliptic functions such as
\begin{equation}
(H_g -\tilde{E}) f_g(x) =0,
\label{eq:Hg}
\end{equation}
where 
\begin{align}
&  H_g= -\frac{d^2}{dx^2} + \sum_{i'=1}^M  \frac{r_{i'} \wp ' (x)}{\wp (x) -\wp (\delta _{i'})} \frac{d}{dx} + \left(l_0 + \sum_{i'=1}^M r_{i'}\right) \left(l_0 +1-  \sum_{i'=1}^M r_{i'}\right) \wp (x) \label{Ino} \\
& \quad \quad +\sum_{i=1}^3 l_i(l_i+1) \wp (x+\omega_i) + \sum_{i'=1}^M \frac{\tilde{s}_{i'}}{\wp (x) -\wp (\delta _{i'})} ,\nonumber \\
& \wp (\delta _{i'})= b_{i'} , \quad ({i'}=1,\dots ,M) .
\end{align}
The parameter $\tilde{s}_{i'}$ $(i'=1,\dots ,M)$ corresponds to the parameter $\tilde{o}_{i'}$, and the parameter $\tilde{p}$ corresponds to $\tilde{E}$.
Apparency of the singularity at $w=\pm \delta _{i'}$ on Eq.(\ref{eq:Hg}) inherits from apparency of the singularity at $w=b_{i'}$ on Eq.(\ref{Feq}).

We now review the propositions on solutions to Eq.(\ref{eq:Hg}) obtained in \cite{TakP}. The first one is an integral representation of solutions in terms of elliptic functions. We set 
$$
\Psi _g(x)=\prod _{i'=1}^M (\wp (x) -\wp (\delta _{i'}))^{r_{i'} /2}.
$$

\begin{prop} \label{prop:integrep} \cite{TakP}
Assume that $l_0, \dots ,l_3 \in \Zint _{\geq 0},$ $r_1, \dots ,  r_k \in \Zint _{\geq 1}$, and the regular singular points $\{ b_1, \dots , b_k \}$ are apparent.
Then there exists an even doubly-periodic function $\Xi (x)$ and a value $Q$ such that
\begin{equation}
\Lambda _g ( x)=\Psi _g (x) \sqrt{\Xi (x)}\exp \int \frac{ \sqrt{-Q}dx}{\Xi (x)}
\label{integ1}
\end{equation}
is a solution to the differential equation (\ref{eq:Hg}).
\end{prop}
For the constructions of $\Xi (x)$ and $Q$, see \cite{TakP}

We now show that a solution to Eq.(\ref{eq:Hg}) can be expressed in the form of the Hermite-Krichever Ansatz. In our situation, the Hermite-Krichever Ansatz asserts that the differential equation has solutions that are expressed as a finite series in the derivatives of an elliptic Baker-Akhiezer function, multiplied by an exponential function. We set
\begin{equation}
\Phi _i(x,\alpha )= \frac{\sigma (x+\omega _i -\alpha ) }{ \sigma (x+\omega _i )} \exp (\zeta( \alpha )x), \quad \quad (i=0,1,2,3),
\label{Phii}
\end{equation}
where $\sigma (x)$ (resp. $\zeta (x)$) is the Weierstrass sigma (resp. zeta) function.
Then we have 
\begin{equation}
\left( \frac{d}{dx} \right) ^{j} \Phi _i(x+2\omega _{k} , \alpha ) = \exp (-2\eta _{k} \alpha +2\omega _{k} \zeta (\alpha )) \left( \frac{d}{dx} \right) ^{j} \Phi _i(x, \alpha )
\label{ddxPhiperiod}
\end{equation}
for $i=0,1,2,3$, $j \in \Zint _{\geq 0}$ and $k=1,3$, where $\eta _k =\zeta (\omega _k)$ $(k=1,3)$.
The following proposition asserts that a solution to Eq.(\ref{eq:Hg}) is written in the form of the Hermite-Krichever Ansatz.
\begin{prop} \label{thm:alpha} \cite{TakP}
(i) Set $\tilde{l} _0 = l_0 +\sum _{i'=1}^M r_{i'}$ and $\tilde{l}_i =l_i$ $(i=1,2,3)$.
The function $\Lambda _g (x)$ in Eq.(\ref{integ1}) is expressed as
\begin{align}
& \Lambda _g (x) = \exp ( \kappa x ) \left( \sum _{i=0}^3 \sum_{j=0}^{\tilde{l}_i-1} \tilde{b} ^{(i)}_j \left( \frac{d}{dx} \right) ^{j} \Phi _i(x, \alpha ) \right)
\label{Lalpha}
\end{align}
for some values $\alpha $, $\kappa$ and $\tilde{b} ^{(i)}_j$ $(i=0,\dots ,3, \: j= 0,\dots ,\tilde{l}_i-1)$, or
\begin{align}
& \Lambda _g (x) =  \exp (\bar{\kappa } x) p(x)
\label{Lalpha0}
\end{align}
for some value $\bar{\kappa }$ and doubly-periodic function $p(x)$. (For a detailed expression of $p(x)$, see \cite{TakP}.)\\
(ii) If $l_0, \dots ,l_3 \in \Zint _{\geq 0},$ $r_1, \dots ,  r_k \in \Zint _{\geq 1}$, and the regular singular points $\{ b_1, \dots , b_k \}$ are apparent, then there exists a non-zero solution to Eq.(\ref{eq:Hg}) that is expressed as Eq.(\ref{Lalpha}) or Eq.(\ref{Lalpha0}).
\end{prop}

The monodromy of the function $\Lambda _g (x)$ is expressed in terms of $\alpha $ and $\kappa $.
In fact, if the function $\Lambda _g (x)$ is written as Eq.(\ref{Lalpha}), then
\begin{align}
& \Lambda _g (x+2\omega _k) = \exp (-2\eta _k \alpha +2\omega _k \zeta (\alpha ) +2 \kappa \omega _k ) \Lambda _g (x) , \quad  (k=1,3). \label{ellint} 
\end{align}

\section{Heun equation} \label{sec:Heun}
For the case $M=0$, Eq.(\ref{Feq}) is regarded as the Heun equation (see Eq.(\ref{Heun})), and it is transformed to the equation
\begin{equation}
\left( -\frac{d^2}{dx^2} + \sum_{i=0}^3 l_i(l_i+1)\wp (x+\omega_i) \right) f(x)= E f(x),\label{Heunell}
\end{equation}
If $l_0, l_1, l_2, l_3 \in \Zint _{\geq 0}$, then the function $\sum_{i=0}^3 l_i(l_i+1)\wp (x+\omega_i)$ is called the Treibich-Verdier potential, and is an example of algebro-geometric finite-gap potential (see \cite{TV,GW,Smi,Tak3}).
For the case $M=0$ and $l_0, l_1, l_2, l_3 \in \Zint _{\geq 0}$, there is no constraint relation for the apparency of additional regular singularity. Hence Propositions \ref{prop:integrep} and  \ref{thm:alpha} holds true.
The function $\Xi (x)$ in Proposition \ref{prop:integrep} is written as 
\begin{equation}
\Xi (x)=c_0(E)+\sum_{i=0}^3 \sum_{j=0}^{l_i-1} b^{(i)}_j (E)\wp (x+\omega_i)^{l_i-j},
\label{XiHeun}
\end{equation}
where the coefficients $c_0(E)$ and $b^{(i)}_j(E)$ are polynomials in $E$, they do not have common divisors and the polynomial $c_0(E)$ is monic. The value $Q$ is expressed as
\begin{align}
 & Q=  \Xi (x)^2\left( E- \sum_{i=0}^3 l_i(l_i+1)\wp (x+\omega_i)\right) +\frac{1}{2}\Xi (x)\frac{d^2\Xi (x)}{dx^2}-\frac{1}{4}\left(\frac{d\Xi (x)}{dx} \right)^2. \label{const}
\end{align}
It follows from Eq.(\ref{Heunell}) that $Q$ is independent of $x$, and it is a monic polynomial in $E$ (see \cite{Tak1}).
A solution to Eq.(\ref{Heunell}) is expressed by an integral (see Eq.(\ref{integ1})), and it is also expressed in a form of the Hermite-Krichever Ansatz (see Proposition \ref{thm:alpha}). It is shown in \cite{Tak4} that the values $\wp (\alpha )$, $\wp ' (\alpha )/\sqrt{-Q}$ and $\kappa /\sqrt{-Q}$ are expressed as rational functions in $E$, and it follows that global monodromy of the Heun equation for the case $l_0, l_1, l_2, l_3 \in \Zint _{\geq 0}$ is written as an elliptic integral.
On the other hand, it is known that global monodromy is also expressed by a hyperelliptic integral (see \cite{Tak3}). By comparing the two expressions, we obtain a hyperelliptic-to-elliptic integral reduction formula (see \cite{Tak4}).

We expressed the functions appeared in this section for the case $l_0=2$, $l_1=l_2=l_3=0$. Note that Eq.(\ref{Heunell}) for the case $l_1=l_2=l_3=0$ is called the Lam\'e equation.

\subsection{The case $M=0$, $l_0=2$, $l_1=l_2=l_3=0$}
The differential equation (see Eq.(\ref{Heunell})) is written as 
\begin{equation}
\left( -\frac{d^2}{dx^2} + 6\wp (x) \right) f(x)= E f(x).\label{Heunell2000}
\end{equation}
Set 
\begin{align}
& \Xi (x)=9\wp(x)^2+3E\wp(x)+E^2-9g_2/4 ,\quad  Q= (E^2-3g_2) \prod _{i=1}^3 (E-3e_i),
\end{align}
where $g_2=-4(e_1e_2+e_2e_3+e_3e_1)$. Then the function
\begin{equation}
\Lambda _g ( x)= \sqrt{\Xi (x)}\exp \int \frac{ \sqrt{-Q}dx}{\Xi (x)},
\end{equation}
is a solution to Eq.(\ref{Heunell2000}).
The monodromy formula in terms of hyperelliptic integral is written as
\begin{equation}
\Lambda _g(x+2\omega _k)=\Lambda _g (x) \exp \left( -\frac{1}{2} \int_{\sqrt{3g_2}}^{E}\frac{ -6\eta _k \tilde{E} +2\omega _k (\tilde{E}^2-3g_2/2) }{\sqrt{-(\tilde{E}^2-3g_2) \prod _{i=1}^3 (\tilde{E}-3e_i)}} d\tilde{E}\right) 
\label{hypellint11000-2}
\end{equation}
for $k=1,3$.

The function $\Lambda _g ( x)$ is also expressed in the form of the Hermite-Krichever Ansatz as
\begin{align}
& \Lambda _g (x) = \exp (\kappa x) \left\{ \bar{b} ^{(0)}_0 \Phi _0 (x, \alpha ) +\bar{b} ^{(0)}_1 \frac{d}{dx} \Phi _0 (x, \alpha ) \right\} 
\end{align}
for the case $E^2\neq 3g_2$, and the values $\alpha $ and $\kappa $ are determined as
\begin{align}
& \wp (\alpha )  = -\frac{E^3-27g_3}{9(E^2-3g_2)}, \quad \kappa = \frac{2}{3} \sqrt{\frac{-(E-3e_1)(E-3e_2)(E-3e_3)}{(E^2-3g_2)}}, \label{kappas42000}
\end{align}
where $g_3=4e_1e_2e_3$. The monodromy is written by using the values $\alpha $ and $\kappa $ (see Eq.(\ref{ellint})). By comparing two expressions of monodromy, we obtain that
\begin{align}
& \int _{\infty} ^{\xi } \frac{d \tilde{\xi }}{\sqrt{4 \tilde{\xi } ^3-g_2 \tilde{\xi } -g_3}} = - \frac{3}{2} \int _{\infty}^{E} \frac{\tilde{E}d\tilde{E}}{\sqrt{-(\tilde{E}^2-3g_2) \prod _{i=1}^3 (\tilde{E}-3e_i)}}, \label{alpint11000^2} \\
&  \kappa = -\frac{1}{2} \int  _{e_i}^{E} \frac{\tilde{E}^2-3g_2/2}{\sqrt{-(\tilde{E}^2-3g_2) \prod _{i=1}^3 (\tilde{E}-3e_i)}}d\tilde{E} + \int _{3e_i} ^{\xi } \frac{ \tilde{\xi } d \tilde{\xi }}{\sqrt{4\tilde{\xi }^3-g_2 \tilde{\xi }-g_3}} ,  \label{kap12311000-2} 
\end{align}
$(i=1,2,3)$ for the transformation 
\begin{equation}
\xi = -\frac{E^3-27g_3}{9(E^2-3g_2)}.
\end{equation}
These formulae reduce hyperelliptic integrals of genus two to elliptic integrals.

\section{Sixth Painlev\'e equation} \label{sec:P6}

We consider the Fuchsian differential equation (\ref{Feq}) for the case $M=1$, $r_1 =1$. Then Eq.(\ref{Feq}) is transformed to 
\begin{equation}
\left\{ -\frac{d^2}{dx^2} + \frac{\wp ' (x)}{\wp (x) -\wp (\delta _{1})} \frac{d}{dx} + \frac{\tilde{s}_{1}}{\wp (x) -\wp (\delta _{1})} +\sum_{i=0}^3 l_i(l_i+1) \wp (x+\omega_i)-\tilde{E}\right\} f_g (x)=0.
\label{HgP6}
\end{equation}
We set
\begin{align}
& b_1 =\wp (\delta _1) , \quad \mu _1= \frac{-\tilde{s}_{1}}{4 b_1^3 -g_2 b_1 -g_3} +\sum _{i=1}^3 \frac{l_i}{2(b_1-e_i)}, \label{pgkr0} \\
& p= \tilde{E} -2(l_1l_2e_3+l_2l_3e_1+l_3l_1e_2) +\sum _{i=1}^3 l_i (l_i e_i +2(e_i+b_1 )) .  \label{pgkr}
\end{align}
The condition that, the regular singular points $x=\pm \delta _1$ is apparent, is written as
\begin{align}
& p=(4 b_1^3 -g_2 b_1 -g_3) \left\{ -\mu _{1} ^2 +\sum_{i=1}^3\frac{l_i+\frac{1}{2}}{b_1-e_i} \mu _{1} \right\} \label{pgkrap} \\
& \quad \quad -b_1 (l_1+l_2+l_3-l_0)(l_1+l_2+l_3+l_0+1) .\nonumber
\end{align}
From now on we assume that $l_0, l_1, l_2, l_3 \in \Zint _{\geq 0}$ and the eigenvalue $\tilde{E}$ satisfies Eqs.(\ref{pgkr}, \ref{pgkrap}). 
Then Propositions \ref{prop:integrep} and \ref{thm:alpha} hold true.
It is known \cite{TakP} that the function $\Xi (x)$ in Proposition \ref{prop:integrep} is written as 
\begin{equation}
\Xi (x)=c_0+\frac{d_0}{(\wp (x)-\wp (\delta _1))}+\sum_{i=0}^3 \sum_{j=0}^{l_i-1} b^{(i)}_j \wp (x+\omega_i)^{l_i-j}.
\label{Fxkr1}
\end{equation}
Ratios of the coefficients $c_0/d_0$ and $b^{(i)}_j/d_0$ $(i=0,1,2,3, \: j=0,\dots ,l_i-1)$ are written as rational functions in variables $b_1$ and $\mu _{1} $.
The value $Q$ in Proposition \ref{prop:integrep} is expressed as a rational function in $b_1$ and $\mu _{1} $ multiplied by $d_0 ^2$. 
By Proposition \ref{thm:alpha}, the eigenfunction $\Lambda _g (x)$ in Eq.(\ref{integ1}) is also expressed in the form of the Hermite-Krichever Ansatz. Namely, it is expressed as 
\begin{align}
& \Lambda _g (x) = \exp \left( \kappa x \right) \left( \sum _{i=0}^3 \sum_{j=0}^{\tilde{l}_i-1} \tilde{b} ^{(i)}_j \left( \frac{d}{dx} \right) ^{j} \Phi _i(x, \alpha ) \right)
\label{LalphaP6}
\end{align}
or
\begin{align}
& \Lambda _g (x) = \exp \left( \bar{\kappa } x \right) p(x)
\label{Lalpha0P6}
\end{align}
for some doubly-periodic function $p(x)$, where $l= l_0+l_1+l_2+l_3 +1$, $\tilde{l} _0 = l_0 +1$ and $\tilde{l}_i =l_i$ $(i=1,2,3)$. For the values $\alpha $ and $\kappa $, we have
\begin{prop} \label{prop:P6HK} \cite{TakP}
Assume that $M=1$, $r_1 =1$, $l_0, l_1, l_2, l_3 \in \Zint _{\geq 0}$ and the value $p$ satisfies Eq.(\ref{pgkrap}). Let $\alpha $ and $\kappa $ be the values determined by the Hermite-Krichever Ansatz  (see Eq.(\ref{LalphaP6})).
Then $\wp (\alpha )$ is expressed as a rational function in variables $b_1$ and $\mu _{1} $, $\wp ' (\alpha )$ is expressed as a product of $\sqrt {-Q}$ and a rational function in variables $b_1$ and $\mu _{1} $, and $\kappa $ is expressed as a product of $\sqrt {-Q}$ and a rational function in variables $b_1$ and $\mu _{1} $.
\end{prop}
If $\alpha \not \equiv 0$ (mod $2\omega _1 \Zint \oplus 2\omega _3 \Zint $), then the function $\Lambda _g (x)$ is expressed as Eq.(\ref{LalphaP6}) and we have
\begin{align}
& \Lambda _g (x+2\omega _k) = \exp (-2\eta _k \alpha +2\omega _j \zeta (\alpha ) +2 \kappa \omega _k ) \Lambda _g (x) , \quad  (k=1,3). \label{ellint00} 
\end{align}

We now discuss the relationship between the monodromy preserving deformation of Fuchsian equations and the sixth Painlev\'e equation. For this purpose we recall other expressions of the sixth Painlev\'e equation. The sixth Painlev\'e equation given as Eq.(\ref{eq:P6eqnintr}) is also written in terms of a Hamiltonian system by adding the variable $\mu$, which is called the sixth Painlev\'e system:
\begin{equation}
\frac{d\lambda }{dt} =\frac{\partial H_{VI}}{\partial \mu}, \quad \quad
\frac{d\mu }{dt} =-\frac{\partial H_{VI}}{\partial \lambda}
\label{eq:Psys}
\end{equation}
with the Hamiltonian 
\begin{align}
H_{VI} = & \frac{1}{t(t-1)} \left\{ \lambda (\lambda -1) (\lambda -t) \mu^2 \right. \label{eq:P6} \\
& \left. -\left\{ \kappa _0 (\lambda -1) (\lambda -t)+\kappa _1 \lambda (\lambda -t) +(\kappa _t -1) \lambda (\lambda -1) \right\} \mu +\kappa (\lambda -t)\right\} ,\nonumber
\end{align}
where $\kappa = ((\kappa _0 +\kappa _1 +\kappa _t -1)^2- \kappa _{\infty} ^2)/4$.
The sixth Painlev\'e equation for $\lambda $ is obtained by eliminating $\mu $ in Eq.(\ref{eq:Psys}).
Set $\omega _1=1/2$, $\omega _3=\tau /2$ and write
\begin{equation}
t= \frac{e_3- e_1}{e_2-e_1}, \quad \lambda = \frac{\wp (\delta )- e_1}{e_2-e_1}.
\end{equation}
Then the sixth Painlev\'e equation is equivalent to the following equation (see \cite{Man,Tks}):
\begin{equation}
\frac{d^2 \delta }{d \tau ^2} = -\frac{1}{4\pi ^2} \left\{ \frac{\kappa _{\infty}^2}{2} \wp ' \left(\delta  \right) + \frac{\kappa _{0}^2}{2} \wp ' \left(\delta +\frac{1}{2} \right) + \frac{\kappa _{1}^2}{2} \wp ' \left(\delta +\frac{\tau +1}{2} \right) +  \frac{\kappa _{t}^2}{2} \wp ' \left(\delta +\frac{\tau }{2}\right) \right\}, \label{eq:P6ellip}
\end{equation}
where $\wp ' (z ) = (\partial /\partial z ) \wp (z)$. In the appendix, we obtain the elliptic form of the sixth Painlev\'e equation (i.e., Eq.(\ref{eq:P6ellip})) from the original sixth Painlev\'e equation (i.e., Eq.(\ref{eq:P6eqnintr})).

It is widely known that the sixth Painlev\'e equation is obtained by the monodnomy preserving deformation of a certain linear differential equation.
Let us introduce the following Fuchsian differential equation:
\begin{equation}
\frac{d^2y}{dw^2} + p_1 (w) \frac{dy}{dw} +p_2(w) y=0, \label{eq:mpdP6}
\end{equation}
where 
\begin{align}
& p_1 (w) = \frac{1-\kappa _0}{w} + \frac{1-\kappa _1}{w-1} + \frac{1-\kappa _t}{w-t} -\frac{1}{w-\lambda}, \\
& p_2 (w) = \frac{\kappa }{w(w-1)} -\frac{t(t-1) H_{VI}}{w(w-1)(w-t)} + \frac{\lambda (\lambda -1) \mu}{w(w-1)(w-\lambda)}.
\end{align}
This equation has five regular singular points $\{ 0,1,t,\infty ,\lambda \}$ and the exponents at $w=\lambda $ are $0$ and $2$.
It follows from Eq.(\ref{eq:P6}) that the regular singular point $w=\lambda $ is apparent.
Then the sixth Painlev\'e equation is obtained by the monodromy preserving deformation of Eq.(\ref{eq:Psys}), i.e., the condition that the monodromy of Eq.(\ref{eq:mpdP6}) is preserved as deforming the variable $t$ is equivalent to that $\mu $ and $\lambda $ satisfy the Painlev\'e system (see Eq.(\ref{eq:Psys})), provided $\kappa _0, \kappa _1, \kappa _t , \kappa _{\infty} \not \in \Zint$. For details, see \cite{IKSY}.

We transform Eq.(\ref{eq:mpdP6}) into the form of Eq.(\ref{HgP6}). We set 
\begin{align}
& w=\frac{\wp (x) -e_1}{e_2-e_1}, \quad y= f_g (x) \prod _{i=1}^3 (\wp (x)-e_i)^{l_i/2}, \label{eq:wwpx} \\
&  \quad t=\frac{e_3-e_1}{e_2-e_1}, \quad \lambda =\frac{b_1 -e_1}{e_2 -e_1}, \quad \wp (\delta _1) =b_1.
\end{align}
Then we obtain Eq.(\ref{HgP6}) by setting
\begin{align}
& \kappa _0 =l_1 +1/2, \quad \kappa _1 =l_2 +1/2, \quad \kappa _t =l_3 +1/2, \quad \kappa _{\infty} =l_0 +1/2, \label{eq:kili} \\
& \mu = (e_2-e_1)\mu _1, \quad \kappa = (l_1+l_2+l_3+l_0 +1)(l_1+l_2+l_3-l_0) , \\
& H_{VI}=\frac{1}{t(1-t)} \left\{ \frac{p+\kappa e_3}{e_2-e_1} +\lambda (1- \lambda)\mu \right\},
\end{align}
(see Eqs.(\ref{pgkr0}--\ref{pgkr})), and Eq.(\ref{eq:P6}) is equivalent to Eq.(\ref{pgkrap}), that means that the apparency of regular singularity is inheritted. Note that the monodromy preserving deformation of Eq.(\ref{eq:mpdP6}) in $t$ corresponds to the monodromy preserving deformation of Eq.(\ref{HgP6}) in $\tau $.

Now we consider the monodromy preserving deformation in the variable $\tau$ ($\omega _1 =1/2, \omega _3=\tau /2$) by applying solutions obtained by the Hermite-Krichever Ansatz for the case $l_i \in \Zint _{\geq 0}$ $(i=0,1,2,3)$.
Let $\alpha $ and $\kappa $ be values determined by the Hermite-Krichever Ansats (see Eq.(\ref{LalphaP6})). We consider the case $Q\neq 0$. Then a basis for solutions to Eq.(\ref{eq:Hg}) is given by $\Lambda _g(x)$ and $\Lambda _g(-x)$, and the monodromy matrix with respect to the cycle $x \rightarrow x+2\omega _k $ ($k=1,3$) is diagonal with the eigenvalues $\exp (\pm ( -2\eta _k \alpha +2\omega _k \zeta (\alpha ) +2 \kappa \omega _k))$ (see Eq.(\ref{ellint00})). 
The values $ -2\eta _k \alpha +2\omega _k \zeta (\alpha ) +2 \kappa \omega _k$ $(k=1,3)$ are preserved by the monodromy preserving deformation.
We set 
\begin{align}
& -2\eta _1 \alpha +2\omega _1 \zeta (\alpha ) +2 \kappa \omega _1 = \pi \sqrt{-1} C_1, \\
& -2\eta _3 \alpha +2\omega _3 \zeta (\alpha ) +2 \kappa \omega _3 = \pi \sqrt{-1} C_3, 
\end{align}
for contants $C_1$ and $C_3$.
By Legendre's relation $\eta _1 \omega _3 - \eta_3 \omega _1 =\pi \sqrt{-1} /2$, we have 
\begin{align}
& \alpha  = C_3 \omega _1 -C_1 \omega _3  \label{al00},\\
& \kappa = \zeta (C_1 \omega _3 -C_3 \omega _1 ) +C_3 \eta _1 -C_1 \eta _3  , \label{kapp00}
\end{align}
From Proposition \ref{prop:P6HK}, the values $\wp (C_3 \omega _1-C_1 \omega _3 ))(=\wp (\alpha))$, $\wp '(C_3 \omega _1-C_1 \omega _3 )/\sqrt {-Q}$ and $(\zeta (C_1 \omega _3 -C_3 \omega _1 ) +C_3 \eta _1 -C_1 \eta _3 )  /\sqrt {-Q}$ are expressed as rational functions in variables $b_1$ and $\mu _1$.
By solving these equations for $b_1$ and $\mu _1$ and evaluating them into Eq.(\ref{HgP6}), the monodromy of the solutions to the differential equation (\ref{HgP6}) on the cycles $x \rightarrow x +2\omega _k$ $(k=1,3)$ are preserved for the fixed values $C_1$ and $C_3$. 
Thus we obtain the following proposition which was established in \cite{TakP}.
\begin{prop} \label{prop:P6} \cite{TakP}
We set $\omega _1= 1/2$, $\omega _3 =\tau /2$ and assume that $l_i \in \Zint _{\geq 0}$ $(i=0,1,2,3)$ and $Q\neq 0$.
By solving the equations in Proposition \ref{prop:P6HK} in variable $b_1 =\wp (\delta _1)$ and $\mu _1$, we express $\wp (\delta _1)$ and $\mu _1$ in terms of $\wp (\alpha)$, $\wp '(\alpha)$ and $\kappa $, and we replace $\wp (\alpha)$, $\wp '(\alpha)$ and $\kappa $ with $\wp (C_3 \omega _1-C_1 \omega _3 )$, $\wp '(C_3 \omega _1-C_1 \omega _3 )$ and $\zeta (C_1 \omega _3 -C_3 \omega _1) +C_3 \eta _1 -C_1 \eta _3$. Then $\delta _1$ satisfies the sixth Painlev\'e equation in the elliptic form
\begin{equation}
\frac{d^2 \delta _1}{d \tau ^2} = -\frac{1}{8\pi ^2} \left\{ \sum _{i=0}^3 (l_i +1/2)^2 \wp '( \delta _1 + \omega _i) \right\}. \label{eq:P6ellipl}
\end{equation}
\end{prop}
We observe the expressions of $b_1$ and $\mu _1$ in detail for the cases $l_0=l_1=l_2=l_3=0$ and $l_0=1$, $l_1=l_2=l_3=0$.

\subsection{The case $M=1$, $r_1 =1$, $l_0=l_1=l_2=l_3=0$}

We investigate the case $M=1$, $r_1 =1$, $l_0=l_1=l_2=l_3=0$ in detail.
The differential equation (\ref{HgP6}) is written as 
\begin{equation}
\left\{ -\frac{d^2}{dx^2} + \frac{\wp ' (x)}{\wp (x) -b_1} \frac{d}{dx} - \frac{ \mu _{1} (4 b_1^3 -g_2 b_1 -g_3)}{\wp (x) -b_1} -p \right\} f_g (x)=0,
\label{Hgkr1l00}
\end{equation}
We assume that $b_1 \neq e_1, e_2, e_3$. The condition that the regular singular points $x= \pm \delta _1$ $(\wp (\delta _{1})=b_1 )$ are apparent is written as
\begin{align}
& p=- (4 b_1^3 -g_2 b_1 -g_3) \mu _{1} ^2 +(6b_1^2 -g_2/2) \mu _{1}  \label{pgkrapl00} 
\end{align}
(see Eq.(\ref{pgkrap})). The doubly-periodic function $\Xi (x)$ (see Eq.(\ref{Fxkr1})) is calculated as 
\begin{equation}
\Xi (x)= 2\mu _1 +\frac{1}{\wp(x)-b_1} .
\end{equation}
The value $Q$ is calculated as
\begin{align}
& Q= 2\mu _1(2\mu _1 (e_1-b_1)+1)(2(e_2-b_1) \mu _1+1)(2\mu _1(e_3-b_1)+1). 
\end{align}
We set
\begin{equation}
\Lambda _g(x) = \sqrt{\Xi (x) (\wp (x) - b_1)} \exp \int \frac{ \sqrt{-Q}dx}{\Xi (x)}.
\label{integ1P6l00}
\end{equation}
%(see Eq.(\ref{integ1P6})).
Then a solution to Eq.(\ref{Hgkr1l00}) is written as $\Lambda _g (x)$, and is expressed in the form of the Hermite-Krichever Ansatz as
\begin{align}
& \Lambda _g (x) = \bar{b} ^{(0)}_0 \exp (\kappa x) \Phi _0 (x, \alpha )
\end{align}
for generic $(\mu_1 , b_1)$.
The values $\alpha $ and $\kappa $ are determined as
\begin{align}
& \wp (\alpha )= b_1 - \frac{1}{2 \mu_1}, \quad \wp '(\alpha )= -\frac{\sqrt{-Q}}{2\mu _1^2} , \quad \kappa = \frac{\sqrt{-Q}}{2\mu _1}.
\end{align}
Hence we have
\begin{align}
& \mu _1 = -\frac{\kappa  }{\wp ' (\alpha )} ,\quad  b_1 = \wp (\alpha ) -\frac{\wp ' (\alpha )}{2\kappa }.
\end{align}
From Proposition \ref{prop:P6}, the function $\delta _1$ determined by
\begin{align}
\wp (\delta _1) =  b_1 & = \wp (C_1 \omega _3 -C_3 \omega_1 ) +\frac{\wp ' (C_1 \omega _3 -C_3 \omega_1 )}{2(\zeta (C_1 \omega _3 -C_3 \omega_1 ) -(C_1 \eta _3 -C_3 \eta _1)) } \label{P6sol0000}
\end{align}
is a solution to the sixth Painlev\'e equation in the elliptic form (see Eq.(\ref{eq:P6ellipl})). This solution coincides with the one found by Hitchin \cite{Hit}.

Now we consider the case $Q=0$. If $Q=0$, then $\mu _1=0$ or $\mu _1 =1/(2(b_1 -e_i))$ for some $i \in \{1,2,3 \}$. 
For the case $\mu _1=0$, the function $\delta _1$, which is determined by
\begin{equation}
\wp (\delta _1) =b_1 =-\frac{D_1 \eta _3 -D_3 \eta _1}{D_1 \omega _3 -D_3 \omega _1},
\label{b1mu10}
\end{equation}
is a solution to the sixth Painlev\'e equation for constants $D_1$ and $D_3$. For the case $\mu _1 =1/(2(b_1 -e_i))$ $(i \in \{1,2,3 \})$, the function $\delta _1$ determined by
\begin{equation}
\wp (\delta _1) =b_1 =\frac{(g_2/4 -2e_i^2)(D_1 \omega _3 -D_3 \omega _1) +e_i (D_1 \eta _3 -D_3 \eta _1) }{e_i (D_1 \omega _3 -D_3 \omega _1) +(D_1 \eta _3 -D_3 \eta _1) }
\label{b1mui}
\end{equation}
is a solution to the sixth Painlev\'e equation.

Eqs.(\ref{b1mu10} ,\ref{b1mui}) are also obtained by suitable limits from Eq.(\ref{P6sol0000}) (see \cite{TakP}), and
the space of the parameters of the solutions to the sixth Painlev\'e equation (i.e. the space of initial conditions) for the case $l_0=l_1=l_2=l_3=0$ is obtained by blowing up four points on the surface $\Cplx /(2\pi \sqrt{-1} \Zint ) \times \Cplx /(2\pi \sqrt{-1} \Zint ) $. This reflects the $A_1 \times A_1 \times A_1 \times A_1$ structure of Riccati solutions by Saito and Terajima \cite{ST}.

\subsection{The case $M=1$, $r_1 =1$, $l_0=1$, $l_1=l_2=l_3=0$}
The differential equation (\ref{HgP6}) for this case is written as 
\begin{equation}
\left\{ -\frac{d^2}{dx^2} + \frac{\wp ' (x)}{\wp (x) -b_1} \frac{d}{dx} - \frac{ \mu _{1} (4 b_1^3 -g_2 b_1 -g_3)}{\wp (x) -b_1} + 2\wp (x) -p \right\} f_g (x)=0,
\label{Hgkr1l01}
\end{equation}
We assume that $b_1 \neq e_1, e_2, e_3$. The condition that the regular singular points $x= \pm \delta _1$ $(\wp (\delta _{1})=b_1 )$ are apparent is written as
\begin{align}
& p=- (4 b_1^3 -g_2 b_1 -g_3) \mu _{1} ^2 +(6b_1^2 -g_2/2) \mu _{1} +2b_1 .\label{pgkrapl01} 
\end{align}
The doubly-periodic function $\Xi (x)$ is calculated as 
\begin{align}
\Xi (x)=& \wp (x) +( (-4b_1^3+b_1g_2+g_3)\mu _1^2+(6b_1^2-g_2/2)\mu _1 -b_1) \\
& +  ((-4b_1^3+b_1g_2+g_3)\mu _1 /2+3b_1^2-g_2/4)/(\wp (x) -b_1) , \nonumber
\end{align}
and the value $Q$ is calculated as
\begin{align}
 Q= -& ((2(4 b_1^3-b_1g_2-g_3)\mu _1^3-(12b_1^2-g_2)\mu _1^2+4) (2(b_1^2+e_1b_1+e_2e_3)\mu _1-2b_1-e_1) \\
& (2(b_1^2+e_2b_1+e_1e_2)\mu _1-2b_1-e_2) (2(b_1^2+e_3b_1+e_1e_3)\mu _1-2b_1-e_3) \nonumber .
\end{align}
We set
\begin{equation}
\Lambda _g(x) = \sqrt{\Xi (x) (\wp (x) - b_1)} \exp \int \frac{ \sqrt{-Q}dx}{\Xi (x)}.
\label{integ1P6l01}
\end{equation}
Then a solution to Eq.(\ref{Hgkr1l01}) is written as $\Lambda _g (x)$, and it is expressed in the form of the Hermite-Krichever Ansatz as
\begin{align}
& \Lambda _g (x) = \exp (\kappa x) \left\{ \bar{b} ^{(0)}_0 \Phi _0 (x, \alpha ) +\bar{b} ^{(0)}_1 \frac{d}{dx} \Phi _0 (x, \alpha ) \right\} 
\end{align}
for generic $(\mu_1 , b_1)$.
The values $\alpha $ and $\kappa $ are determined as
\begin{align}
& \wp (\alpha )  = \frac{2(4b_1^3-b_1g_2-g_3)b_1 \mu _1^3+(-24b_1^3+4g_2b_1+3g_3)\mu _1^2+(24b_1^2-2g_2)\mu _1-8b_1}{2(4b_1^3-b_1 g_2-g_3)\mu _1^3 -(12b_1^2-g_2)\mu _1^2+4},\\
& \wp '(\alpha )= \frac{-4((4b_1^3-b_1g_2-g_3)\mu _1^3-(12b_1^2-g_2)\mu _1^2+12b_1\mu _1-4)}{(2(4b_1^3-b_1 g_2-g_3)\mu _1^3 -(12b_1^2-g_2)\mu _1^2+4)^2}\sqrt{-Q}  ,\nonumber \\
& \kappa = \frac{2\mu _1}{2(4b_1^3-b_1 g_2-g_3)\mu _1^3 -(12b_1^2-g_2)\mu _1^2+4}\sqrt{-Q}.\nonumber 
\end{align}
Hence we have
\begin{align}
& b_1 = \frac{2\wp (\alpha ) \kappa ^3-3\wp '(\alpha )\kappa ^2+(6\wp (\alpha ) ^2 -g_2)\kappa -\wp (\alpha ) \wp '(\alpha )}{2(\kappa ^3-3\wp (\alpha ) \kappa +\wp '(\alpha ))} ,\\
& \mu _1 = \frac{2(\kappa ^3-3\wp (\alpha ) \kappa +\wp '(\alpha ))\kappa }{-2\wp '(\alpha )\kappa ^3+(12\wp (\alpha ) ^2-g_2)\kappa ^2-6\wp (\alpha ) \wp '(\alpha )\kappa +\wp '(\alpha )^2}.
\end{align}
From Proposition \ref{prop:P6}, the function $\delta _1$ determined by
\begin{align}
& \wp (\delta _1) = b_1 = \label{P6sol1000} \\ 
& \quad  \frac{2\wp (\omega ) (\zeta (\omega )- \eta )^3+3\wp '(\omega )(\zeta (\omega )- \eta )^2+(6\wp (\omega ) ^2 -g_2)(\zeta (\omega )- \eta )+\wp (\omega ) \wp '(\omega )}{2((\zeta (\omega )- \eta )^3-3\wp (\omega ) (\zeta (\omega )- \eta ) -\wp '(\omega ))} , \nonumber \\
& (\omega = C_1 \omega _3 -C_3 \omega_1 , \quad  \eta = C_1 \eta _3 -C_3 \eta _1) ,\nonumber
\end{align}
is a solution to the sixth Painlev\'e equation in the elliptic form (see Eq.(\ref{eq:P6ellipl})).
In the sixth Painlev\'e equation, it is known that the case $(\kappa _{0}, \kappa _{1}, \kappa _{t}, \kappa _{\infty}) =(1/2, 1/2, 1/2, 3/2 ) $ is linked to the case $(\kappa _{0}, \kappa _{1}, \kappa _{t}, \kappa _{\infty}) =(1/2, 1/2, 1/2, 1/2 ) $ by B\"acklund transformation (see \cite{TOS}).
By transformating the solution in Eq.(\ref{P6sol0000}) of the case $(\kappa _{0}, \kappa _{1}, \kappa _{t}, \kappa _{\infty}) =(1/2, 1/2, 1/2, 1/2 )$ to the one of the case  $(\kappa _{0}, \kappa _{1}, \kappa _{t}, \kappa _{\infty}) =(1/2, 1/2, 1/2, 3/2 )$, we recover the solution in Eq.(\ref{P6sol1000}).

Now we consider the case $Q=0$. If $Q=0$, then $\mu_1  $ is a solution to the equation $2(4 b_1^3-b_1g_2-g_3)\mu _1^3-(12b_1^2-g_2)\mu _1^2+4 =0$ or $\mu _1=(2b_1+e_i)/(2(b_1^2+e_ib_1+e_i^2-g_2/4))$ for some $i \in \{1,2,3 \}$.
We set $\omega = D_1 \omega _3 -D_3 \omega_1$ and $\eta = D_1 \eta _3 -D_3 \eta _1$, where $D_1$ and $D_3$ are constants. For the case that $\mu_1  $ is a solution to the equation $2(4 b_1^3-b_1g_2-g_3)\mu _1^3-(12b_1^2-g_2)\mu _1^2+4 =0$, the corresponding solutions to the sixth Painlev\'e equation are written as the function $\delta _1$, where 
\begin{equation}
\wp( \delta _1)= b_1 = \frac{4 \eta ^3 +g_2 \omega ^2 \eta -2 g_3 \omega ^3}{\omega (g_2 \omega ^2-12\eta ^2 )}. \label{b1mu1l01}
\end{equation}
For the case $\mu _1=(2b_1+e_i)/(2(b_1^2+e_ib_1+e_i^2-g_2/4))$ ($i \in \{1,2,3 \}$), we have 
\begin{equation}
\wp( \delta _1)= b_1 = \frac{-g_2 e_i \omega /2 +(6e_i^2 -g_2)\eta}{(6e_i^2 -g_2)\omega -6e_i\eta }. \label{b1muil01}
\end{equation}

\section{Summary and concluding remarks} \label{sec:rmk}

The Heun equation (see Eq.(\ref{Heun})) is the standard canonical form of a Fuchsian equation with four singularities.
By transforming the Heun equation to the form of elliptic functions, we find that solving the Heun equation is equivalent to investigating spectral and eigenstates of quantum $BC_1$ Inozemtsev model (see Eq.(\ref{Heunell})). Note that the Hamiltonian of the quantum $BC_1$ Inozemtsev model \cite{Ino} is given by
\begin{equation}
H= -\frac{d^2}{dx^2} + \sum_{i=0}^3 l_i(l_i+1)\wp (x+\omega_i) .
\end{equation}

On the other hand, by adding an apparent singularity to the Heun equation, we obtain Fuchsian differential equations that produce the sixth Painlev\'e equation by monodromy preserving deformation (see section \ref{sec:P6}).

The sixth Painlev\'e equation is rewritten as Eq.(\ref{eq:P6ellip}), and it is equivalent to the Hamiltonian system (see \cite{Tks})
\begin{align}
& 2\pi \sqrt{-1} \frac{d \delta }{d\tau } = \frac{\partial {\mathcal H}}{\partial \gamma}, \quad \quad 2\pi \sqrt{-1} \frac{d \gamma}{d\tau } = -\frac{\partial {\mathcal H}}{\partial \delta }, \label{eq:P6ellHsys} \\
& {\mathcal H}= \frac{1}{2} \left( \gamma ^2- \sum _{i=0}^3 (l_i +1/2)^2 \wp ( \delta + \omega _i) \right).
\end{align}
If we replace $2\pi \sqrt{-1} \frac{d }{d\tau } $ by $ \frac{d }{ds }$ ($s$: time variable, independent of $\tau$) formally, we obtain the classical $BC_1$ Inozemtsev system \cite{Ino}. In other words, the sixth Painlev\'e equation is a non-autonomous version of the classical $BC_1$ Inozemtsev system.
To summarize, we present the following diagram.

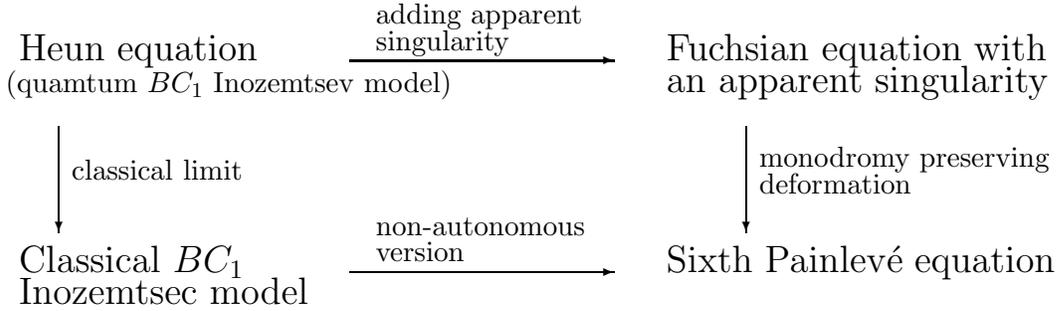
\begin{picture}(370,130)
\put(5,100){\large{Heun equation}} 
\put(0,88){\small{(quamtum $BC_1$ Inozemtsev model)}}

\put(250,100){\large{Fuchsian equation with}}
\put(250,88){\large{an apparent singularity}}

\put(5,20){\large{Classical $BC_1$}}
\put(5,8){\large{Inozemtsec model}}

\put(250,20){\large{Sixth Painlev\'e equation}}

\put(130,100){\vector(1,0){100}}
\put(140,114){\small{adding apparent}}
\put(140,104){\small{singularity}}

\put(130,20){\vector(1,0){100}}
\put(140,34){\small{non-autonomous}}
\put(140,24){\small{version}}

\put(20,75){\vector(0,-1){40}}
\put(25,55){\small{classical limit}}

\put(280,75){\vector(0,-1){40}}
\put(285,60){\small{monodromy preserving}}
\put(285,50){\small{deformation}}

\end{picture}

Before starting this work, the author noticed that the parameters, that the monodromy of solutions to the Heun equation have expressions in terms of elliptic or hyperelliptic integrals, resemble the ones in the sixth Painlev\'e equation that has explicit two-parameter solutions.
Typical two-parameter solutions are Picard's and Hitchin's solutions.
In this paper, we partially obtain an explanation of this phenomena by intermediating Fuchsian differential equations with an apparent singularity, though the corresponding parameters on the sixth Painlev\'e equation are a little off as $O_1$ and $O_1 \cup O_2$, where
\begin{align}
& O_1= \left\{ (\kappa _0 , \kappa _1 , \kappa _t, \kappa _{\infty} ) | 
\kappa _0 , \kappa _1 , \kappa _t, \kappa _{\infty} \in \Zint +\frac{1}{2} \right\}, \\
& O_2 = \left \{(\kappa _0 , \kappa _1 , \kappa _t, \kappa _{\infty} ) \left| 
\begin{array}{ll}
\kappa _0 , \kappa _1 , \kappa _t, \kappa _{\infty} \in \Zint \\
\kappa _0 + \kappa _1 + \kappa _t + \kappa _{\infty}  \in 2 \Zint 
\end{array}
\right. \right\}. 
\end{align}
For the case $(\kappa _0 , \kappa _1 , \kappa _t, \kappa _{\infty} ) \in O_1$, solutions of the linear differential equation are investigated by our method, and solutions of the sixth Painlev\'e equation follow from them (see Proposition \ref{prop:P6}).
By B\"acklund transformation of the sixth Painlev\'e equation (see \cite{TOS} etc.), Hitchin's solution (i.e., solutions for the case $(\kappa _0 , \kappa _1 , \kappa _t, \kappa _{\infty} )=(\frac{1}{2}, \frac{1}{2}, \frac{1}{2}, \frac{1}{2})$) is transformed to the solutions for the case $(\kappa _0 , \kappa _1 , \kappa _t, \kappa _{\infty} ) \in O_1 \cup O_2$.
But we cannot obtain results on integral representation and the Hermite-Krichever Ansatz by our method for the case $(\kappa _0 , \kappa _1 , \kappa _t, \kappa _{\infty} ) \in O_2$.
Note that the condition $(\kappa _0 , \kappa _1 , \kappa _t, \kappa _{\infty} ) \in O_2$ corresponds to the condition $l_0, \dots ,l_3 \in \Zint +\frac{1}{2}$, $l_0+ l_1 +l_2 + l_3 \in 2 \Zint $.
How can we investigate solutions and their monodromy of the linear differential equation for the cases $l_0, \dots ,l_3 \in \Zint +\frac{1}{2}$, $l_0+ l_1 +l_2 + l_3 \in 2 \Zint $?

\appendix
\section{Elliptic form of sixth Painv\'e equation} \label{sect:append}
We calculate the differentiation of modular functions, which will be used to rewrite the sixth Painlev\'e equation.
\begin{prop} (c.f. \cite{Man}) \label{prop:mod1}
Set $\omega _1 =1/2$, $\omega _3=\tau /2$, $t=(e_3-e_1)/(e_2-e_1)$. Then we have 
\begin{align}
& \frac{dt}{d\tau } =\frac{(e_2-e_1)t(t-1)}{\pi \sqrt{-1}}, \label{eq:dtdtau} \\
& \frac{d}{d\tau } \left( \frac{1}{(e_2-e_1)^{1/2}} \right) = \frac{\eta _1 +e_3 /2}{\pi \sqrt{-1} (e_2-e_1)^{1/2}}, \\
& \frac{d}{d\tau } \left((e_2-e_1)^{\alpha } \right) = - \frac{\alpha (2\eta _1 +e_3 )(e_2-e_1)^{\alpha }}{\pi \sqrt{-1}} .
\end{align}
\end{prop}
\begin{proof}
Set $u=(z-e_1)/(e_2-e_1)$. Then 
\begin{align}
& \frac{\tau}{2} = \int _{1/2}^{(1+\tau )/2} dx= \int _{e_1} ^{e_2} \frac{dz}{\wp '(x)} \label{eq:inte1e2} \\
& = \int _{e_1} ^{e_2} \frac{dz}{2\sqrt{(z-e_1)(z-e_2)(z-e_3)}} = \frac{1}{2(e_2-e_1)^{1/2}} \int _0^1 \frac{du}{\sqrt{u(u-1)(u-t)}}, \nonumber \\
& -\frac{1}{2} = \int _{e_1} ^{\infty } \frac{dz}{2\sqrt{(z-e_1)(z-e_2)(z-e_3)}} = \frac{1}{2(e_2-e_1)^{1/2}} \int _0^{\infty} \frac{du}{\sqrt{u(u-1)(u-t)}}, \label{eq:inte1inf}
\end{align}
By differentiating Eq.(\ref{eq:inte1e2}) in variable $\tau $, we have
\begin{align}
\frac{1}{2} = & \frac{d}{d\tau } \left( \frac{1}{2(e_2-e_1)^{1/2}} \right) \int _0^1 \frac{du}{\sqrt{u(u-1)(u-t)}} \\
& + \frac{1}{2} \left( \frac{dt}{d\tau } \right) \frac{1}{2(e_2-e_1)^{1/2}} \int _0^1 \frac{du}{(u-t)\sqrt{u(u-1)(u-t)}} , \nonumber
\end{align}
and it follows from Eq.(\ref{eq:inte1e2}) that
\begin{align}
& \frac{1}{2(e_2-e_1)^{1/2}} \int _0^1 \frac{du}{(u-t)\sqrt{u(u-1)(u-t)}}  = (e_2 -e_1) \int _{e_1} ^{e_2} \frac{dz}{(\wp (x)- e_3) \wp '(x)} \label{eq:int01} \\
& = \frac{e_2-e_1}{(e_3-e_2)(e_3-e_1 )} \int _{1/2}^{(1+\tau )/2} (\wp (x + \tau/2) -e_3) dx \nonumber \\
& = \frac{e_2-e_1}{(e_3-e_2)(e_3-e_1 )} (-\eta (1/2+\tau) + \eta (1/2 +\tau /2) -e_3 \tau /2 )  \nonumber \\
& = \frac{e_2-e_1}{(e_3-e_2)(e_3-e_1 )} (-\eta _3 -e_3 \tau /2 ) . \nonumber
\end{align}
Hence 
\begin{align}
& \frac{1}{2} = (e_2-e_1)^{1/2} \tau  \frac{d}{d\tau } \left( \frac{1}{2(e_2-e_1)^{1/2}} \right)  + \frac{1}{2} \left( \frac{dt}{d\tau } \right) \frac{(e_2-e_1)(-\eta _3 -e_3 \tau /2) }{(e_3-e_2)(e_3-e_1 )} . 
\end{align}
Similarly it follows from differentiating Eq.(\ref{eq:inte1inf}) that
\begin{align}
& 0 = - (e_2-e_1)^{1/2} \frac{d}{d\tau } \left( \frac{1}{2(e_2-e_1)^{1/2}} \right)  + \frac{1}{2} \left( \frac{dt}{d\tau } \right) \frac{(e_2-e_1)(\eta _1 +e_3 /2) }{(e_3-e_2)(e_3-e_1 )} . \label{eq:dtaue2e1}
\end{align}
From these equalities we have
\begin{equation}
\left( \frac{dt}{d\tau } \right) \frac{(e_2-e_1)(\tau \eta _1 -\eta _3 ) }{(e_3-e_2)(e_3-e_1 )} =1 .
\end{equation}
By Legendre's relation $\eta _1 \tau -\eta _3 = \pi \sqrt{-1}$ and definition of $t$, it follows that
\begin{equation}
\frac{dt}{d\tau } =\frac{(e_2-e_1)t(t-1)}{\pi \sqrt{-1}}.
\end{equation}
Combining with Eq.(\ref{eq:dtaue2e1}), we obtain that 
\begin{equation}
\frac{d}{d\tau } \left( \frac{1}{(e_2-e_1)^{1/2}} \right) = \frac{\eta _1 +e_3 /2}{\pi \sqrt{-1} (e_2-e_1)^{1/2}}.
\end{equation}
The derivation of the function $(e_2-e_1)^{\alpha }$ is calculated as
\begin{equation}
\frac{d}{d\tau } (e_2-e_1)^{\alpha }= -2\alpha (e_2-e_1)^{\alpha +1/2} \frac{d}{d\tau } \left( \frac{1}{(e_2-e_1)^{1/2}} \right) = - \frac{\alpha (2\eta _1 +e_3 )(e_2-e_1)^{\alpha }}{\pi \sqrt{-1}} .
\end{equation}
\end{proof}
\begin{prop}
Set $g_2=-4(e_1e_2+e_2e_3+e_3e_1)$. We have
\begin{align}
& \frac{de_i}{d\tau} = \frac{-2\eta _1 e_i +e_i^2-g_2/6}{\pi \sqrt{-1}} , \label{eq:deidtau}
\end{align}
for $i=1,2,3$.
\end{prop}
\begin{proof}
It follows from Proposition \ref{prop:mod1} that
\begin{equation}
\frac{d}{d\tau } (e_2-e_1)^{\pm 1}= \mp  \frac{(2\eta _1 +e_3 )(e_2-e_1)^{\pm 1}}{\pi \sqrt{-1}} . \label{eq:e21}
\end{equation}
Combining with Eq.(\ref{eq:dtdtau}), we obtain that
\begin{equation}
\frac{d}{d\tau } (e_3-e_1)= - \frac{(2\eta _1 +e_2 )(e_3-e_1)}{\pi \sqrt{-1}} . \label{eq:e31}
\end{equation}
By adding two functions, we have 
\begin{equation}
\frac{d}{d\tau } (-3e_1)= - \frac{-6e_1 \eta _1 +2e_2e_3 -e_1^2}{\pi \sqrt{-1}} .\label{eq:-3e1}
\end{equation}
Hence, we obtain Eq.(\ref{eq:deidtau}) for the case $i=1$. 
Eq.(\ref{eq:deidtau}) for the case $i=2$ (resp. $i=3$) follows from Eqs.(\ref{eq:-3e1}, \ref{eq:e21}) (resp. Eqs.(\ref{eq:-3e1}, \ref{eq:e31})).
\end{proof}

\begin{prop}
\begin{equation}
\frac{d\eta _1}{d\tau } = \frac{-\eta _1^2 +g_2/48}{\pi \sqrt{-1}}. \label{eq:detadtau}
\end{equation}
\end{prop}
\begin{proof}
It follows similarly to Eq.(\ref{eq:int01}) that
\begin{equation}
\frac{1}{2} \int _0 ^{\infty } \frac{du}{(u-t)\sqrt{u(u-1)(u-t)}} = \frac{(e_2-e_1)^{3/2}}{(e_3-e_2)(e_3-e_1) }(\eta _1+e_3/2). \label{eq:int0infty}
\end{equation}
We differentiate Eq.(\ref{eq:int0infty}) in $t$. From the l.h.s, we have
\begin{align*}
& \frac{3}{2} \int _0 ^{\infty } \frac{du}{2(u-t)^2\sqrt{u(u-1)(u-t)}} = \frac{3}{2} \frac{(e_2-e_1)^{5/2}}{(e_3-e_2)^2(e_3-e_1)^2 } \int _{1/2} ^0 (\wp (x+\tau/2 )-e_3)^2 dx , \\
& \int _{1/2} ^0 (\wp (x+\tau/2 )-e_3)^2dx= \int _{1/2} ^0 \left( \frac{\wp '' (x+\tau/2 )^2}{6} -2 e_3 \wp (x+\tau/2 ) + e_3^2 +\frac{g_2}{12} \right) dx \\
& = -2e_3\eta _1 -\frac{1}{2}\left( e_3^2+\frac{g_2}{12} \right) .\\
\end{align*}
From the r.h.s, we have
\begin{align*}
& \frac{(e_2-e_1)^{3/2}}{(e_3-e_2)^2(e_3-e_1)^2 }\left\{ \left( \eta _1 +\frac{e_3}{2} \right) \left( \eta _1 -\frac{5e_3}{2} \right) +\pi \sqrt{-1} \frac{d}{d\tau} \left( \eta _1 +\frac{e_3}{2} \right) \right\}.
\end{align*}
Hence, we obtain 
\begin{equation}
\frac{d}{d\tau } \left( \eta _1 +\frac{e_3}{2} \right)= \frac{1}{\pi \sqrt{-1}} \left\{ -\eta ^2 -\eta _1 e_3 + \frac{e_3^2}{2} -\frac{g_2}{16} \right\} ,
\end{equation}
and Eq.(\ref{eq:detadtau}).
\end{proof}

We now we show that the sixth Painlev\'e equation (see Eq.(\ref{eq:P6eqnintr})) can be rewritten to an elliptic form (see Eq.(\ref{eq:P6ellip})).

\begin{prop} \cite{Man}
Set
\begin{equation}
\omega _1 =1/2, \quad \omega _3 =\tau/2, \quad t=\frac{e_3-e_1}{e_2-e_1}, \quad \lambda =\frac{\wp (\delta ) -e_1}{e_2-e_1}.
\end{equation}
Then the sixth Painlev\'e equation 
\begin{align}
\frac{d^2\lambda }{dt^2} = & \frac{1}{2} \left( \frac{1}{\lambda }+\frac{1}{\lambda -1}+\frac{1}{\lambda -t} \right) \left( \frac{d\lambda }{dt} \right) ^2 -\left( \frac {1}{t} +\frac {1}{t-1} +\frac {1}{\lambda -t} \right)\frac{d\lambda }{dt} \label{eq:pP6} \\
& +\frac{\lambda (\lambda -1)(\lambda -t)}{t^2(t-1)^2}\left\{ \frac{\kappa _{\infty}^2}{2} -\frac{\kappa _{0}^2}{2}\frac{t}{\lambda ^2} +\frac{\kappa _{1}^2}{2}\frac{(t-1)}{(\lambda -1)^2} +\frac{(1-\kappa _{t}^2)}{2}\frac{t(t-1)}{(\lambda -t)^2} \right\} \nonumber
\end{align}
is equivalent to the equation
\begin{equation}
\frac{d^2 \delta }{d \tau ^2} = -\frac{1}{4\pi ^2} \left\{ \frac{\kappa _{\infty}^2}{2} \wp ' \left(\delta  \right) + \frac{\kappa _{0}^2}{2} \wp ' \left(\delta +\frac{1}{2} \right) + \frac{\kappa _{1}^2}{2} \wp ' \left(\delta +\frac{\tau +1}{2} \right) +  \frac{\kappa _{t}^2}{2} \wp ' \left(\delta +\frac{\tau }{2}\right) \right\}. \label{eq:pP6el}
\end{equation}
\end{prop}
\begin{proof}
It follows from the relation $\lambda =(\wp (\delta ) -e_1)/(e_2-e_1)$ that
\begin{equation}
\delta = \int_0 ^\delta dx = \int _{\infty} ^{\lambda } \frac{e_2-e_1}{\wp ' (x)} du = \frac{1}{2(e_2-e_1)^{1/2}} \int _{\infty} ^{\lambda } \frac{du}{\sqrt{u(u-1)(u-t)}} . \label{eq:deltaint}
\end{equation}
We differentiate Eq.(\ref{eq:deltaint}) by the variable $\tau$. Then we have
\begin{align}
& \frac{d\delta }{d\tau } = \frac{\eta_1 +e_3/2}{\pi \sqrt{-1}}\delta + \frac{(e_2-e_1)^{1/2}t(t-1)}{2\pi \sqrt{-1}}\left\{ \frac{d\lambda }{dt}\frac{1}{\sqrt{\lambda (\lambda -1)(\lambda -t)}} \right. \label{eq:deltaint1} \\
& \left. \quad \quad \quad \quad \quad  \quad \quad \quad \quad \quad  \quad \quad \quad \quad \quad +\frac{1}{2} \int _{\infty} ^{\lambda } \frac{du}{(u-t)\sqrt{u(u-1)(u-t)}} \right\} . \nonumber
\end{align}
Note that we used Proposition \ref{prop:mod1}. We differentiate Eq.(\ref{eq:deltaint1}) once more. By applying formulae on the differentiation of modular functions, we obtain that
\begin{align}
& \frac{d^2\delta }{d\tau ^2} = \frac{t^2(t-1)^2(e_1-e_3)^{3/2}}{-2\pi ^2} \left[ \frac{1}{\sqrt{\lambda (\lambda -1)(\lambda -t)}}\left\{ \frac{d^2\lambda }{dt ^2} \right. \right. \label{eq:d2ddtau} \\
& \quad \quad \left.  -\frac{1}{2}\left(\frac{1}{\lambda }+\frac{1}{\lambda -1}+\frac{1}{\lambda -t}\right) \left( \frac{d\lambda }{dt} \right) ^2+ \left(\frac{1}{t }+\frac{1}{t -1}+\frac{1}{\lambda -t}\right)\frac{d\lambda }{dt} \right\} \nonumber \\
& \quad \quad \left. + \frac{1}{4t(t-1)} \int _{\infty} ^{\lambda } \frac{(u^2+2tu-2u-t) du}{(u-t)^2\sqrt{u(u-1)(u-t)}} \right] ,\nonumber
\end{align}
and we have
\begin{equation}
\int _{\infty} ^{\lambda } \frac{(u^2+2tu-2u-t) du}{(u-t)^2\sqrt{u(u-1)(u-t)}} = -2 \sqrt{\frac{\lambda (\lambda -1)}{(\lambda -t)^3}}.
\end{equation}
It follows from $\lambda =(\wp (\delta ) -e_1)/(e_2-e_1)$ that
\begin{align}
& \wp ' (\delta )=2 (e_1-e_2)^{3/2} \sqrt{\lambda (\lambda -1)(\lambda -t)}, \quad \frac{\wp ' (\delta +1/2)}{\wp ' (\delta )} =-\frac{t}{\lambda ^2} , \label{eq:wpd} \\
& \frac{\wp ' (\delta +(\tau +1)/2)}{\wp ' (\delta )} =\frac{t-1}{(\lambda -1)^2} , \quad \frac{\wp ' (\delta +\tau /2)}{\wp ' (\delta )} =\frac{t(1-t)}{(\lambda -t)^2} .\nonumber
\end{align}
By combining Eqs.(\ref{eq:d2ddtau}-\ref{eq:wpd}), we obtain the equivalence of Eq.(\ref{eq:pP6}) and Eq.(\ref{eq:pP6el}).
\end{proof}

\end{document}